# The Effect of Using Concept Maps on Advancing Students' Conceptual Understanding of Euler Circuit

Derar Serhan[#1], Muhammed Syam [*2], Qasem AlMdallal[#3]

[#1] *Mathematics Education & Information Science, ECAE, UAE*
[#2] *Dept. of Mathematical Sciences, College of Science, UAEU, UAE*
[#3] *Dept. of Mathematical Sciences, College of Science, UAEU, UAE*

*Abstract*— This study investigated the use of concept maps on advancing students' conceptual understanding of the Euler Circuit concept. The participants were students at the college of Humanities and Social Sciences with a background in Art. A pre-treatment test was given to students during the eleventh week of the semester, while a post-treatment test was given during the fourteenth week. The study reports the findings from the responses of eighty three participants who took both the pre- and post-treatment test. The results of this study indicated that students, after using concept maps, had a richer conceptual understanding of the Euler Circuit concept and were capable of constructing more representations of this concept.

*Keywords*— **Euler Circuit, Concept Maps.**

## I. INTRODUCTION

The impetus for this study was born when students' responses to the final exam in a previous semester showed that they did not properly understand the Euler Circuit concept. For example, students were required to respond to the following question in their final exam:

- Which of the following graphs has an Euler Circuit?

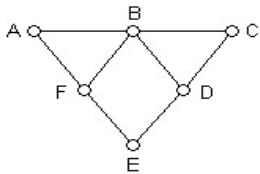 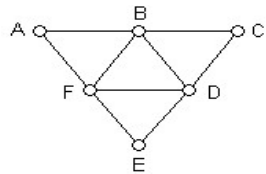 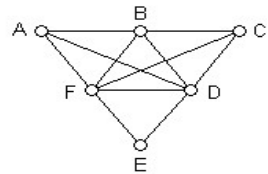

Graph 1    Graph 2    Graph 3

- Graph 2 only
- Graph 1 only
- Graphs 1 and 3
- Graph 3 only

Most of the students (more than 80%) did not answer this question correctly. Because of that, we looked for methods to improve students' conceptual understanding of the Euler Circuit. Since concept maps allow for making connections between all concepts that relate to the Euler Circuit, they were chosen as a treatment tool in this study. The purpose of this study was to investigate the use of concept maps on advancing students' understanding of the Euler Circuit concept.

## II. MATHEMATICAL BACKGROUND

In the 1700s, a famous story (*the Königsberg bridges story*) started in the city of Kaliningrad, then called Königsberg. Through the city ran a river, The Pregel River, which divided the city into two islands A and B, with seven bridges connecting the north and south banks of the river (R and L respectively). Young people at that time used to play a game, its aim was to take a walk around town fully crossing every bridge once and only once. They tried many times but were unable to do so. Leonhard Euler (a famous mathematician) studied that problem and proved that it was impossible to take a walk around town fully crossing every bridge once and only once. Euler redrew the map of the bridges as follows:

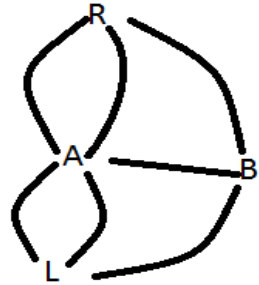

Fig. 1 Mathematical graph of the Pregel River

In this graph, Euler defined the vertices to be R, A, B, L and the edges to be the seven bridges. Two edges are adjacent if they have a common vertex. For example, the edge AB and the edge BL are adjacent since B is a common vertex between them. The graph in Fig.1 is connected, i.e.; you can travel from any vertex to another by a sequence of adjacent edges. A path is a sequence of vertices such that each vertex in the sequence is adjacent to the next one. An Euler path is a path that passes through every edge of a graph once and only once. A Circuit is a sequence of vertices such that each vertex in the sequence is adjacent to the next one and it starts and ends at the same vertex. An Euler Circuit is a Circuit that passes through every edge of a graph.

In this case, we can look at the Königsberg bridges puzzle as if there is an Euler path to the graph of Fig.1. If we want to start and end at the same vertex, then an Euler Circuit is applied to the graph of Fig.1. Euler classified the vertices into odd and even based on the number of edges at that vertex. If the number of edges is odd, then it is called odd vertex. Otherwise, it is called even vertex. Based on that, Euler stated the following theorem:

**Theorem 1:**
- If a graph is *connected*, and every vertex is *even*, then it has an Euler Circuit (at least one, usually more).
- If a graph has *any* odd vertices, then it does not have an Euler Circuit.
- If a graph is *connected*, and has exactly *two* odd vertices, then it has an Euler path. Any such path must start at one of the odd vertices and end at the other one.
- If a graph has *more than two* odd vertices, then it cannot have an Euler path.

The graph in Fig.1 is connected and degrees of A, B, L, and R are 5, 3, 3, and 3, respectively. From Theorem 1, the graph in Fig.1 neither has Euler Circuit nor Euler path. This means Königsberg bridges puzzle does not have a solution.

III. LITERATURE REVIEW

Since the purpose of this study was to investigate the use of concept maps on advancing students' conceptual understanding, it is important to clarify the meaning of understanding and the use of concept maps.

Reference [1] structured their discussion about the meaning of knowing and understanding mathematics around the following themes: understanding as representation; understanding as knowledge structures; understanding as connections among types of knowledge; and understanding as situated cognition.

Reference [2] emphasized understanding as connections between conceptual and procedural knowledge. They defined conceptual knowledge as "knowledge that is rich in relationships. It can be thought of as a connected web of knowledge, a network in which the linking relationships are as prominent as the discrete pieces of information" (pp. 3-4). On the other hand, they defined procedural knowledge as knowledge that is "made up of two distinct parts. One part is composed of the formal language, or symbol representation system, of mathematics. The other part consists of the algorithms, or rules, for completing mathematical tasks"(p. 6). For students to develop a strong mathematical knowledge, it is very important for them to link conceptual and procedural knowledge, without developing such a link between the two types of knowledge, students will not be fully proficient in dealing with mathematical concepts.

Concept maps are used for representing knowledge relationships, they were used for different purposes and in different disciplines. Concept maps were used as an instructional tool [3]-[4], and also as a tool to assess students' understanding of a given concept [5]. Concept maps has been used in many disciplines to assess student learning, it has been used to assess students' learning and understanding in psychology [6], in statistics [7], and in mathematics [8].

In addition to that, concept mapping has also been used to teach classroom content, [9] studied the effectiveness of using concept maps as an instructional tool in a Calculus class. They used concept maps to enhance students' understanding of the concept of the derivative at a point. Moreover, concept maps were used in teaching science and biology [10]-[13].

A concept map is a visual representation of a focal concept that is connected with other related concepts. To develop a concept map, key concepts and the relationships among them are identified and arranged around a focal concept where lines connect between the focal concept and the key concepts around it. These lines are labelled in ways that describe the relationships between them [14].

IV. PURPOSE OF THE STUDY

The purpose of this study was to investigate the effect of using concept maps on advancing students' understanding of the Euler Circuit concept. Euler Circuit is an important concept in graph theory. It depends on many basic concepts such as vertex, edge, path, and Circuit. Students use it to solve many real life problems. One of these problems is Königsberg bridges puzzle. Students face difficulty understanding this concept and they are unable to make connections between its different representations. The use of concept maps can help students make these connections and organize their own thoughts in a graphical way that is easy to follow. The present study compared between the students' conceptual understanding of the Euler Circuit before and after the use of concept maps. The main research questions was the following:

> Is there a difference between students' conceptual understanding of the Euler Circuit before and after the use of concept maps?

V. METHOD

A. *Participants*

The students who participated in this study were enrolled in a Contemporary of Applied Mathematics course at a major university in the United Arab Emirates. The participants, who had an Art background, were students at the college of Humanities and Social Sciences. One hundred students took the pre-treatment test and only eighty three of them took the post-treatment test. For this reason, we excluded the responses of seventeen students who did not take the post-treatment test. The students in this course used the book "*Excursions in Modern Mathematics*" [15]. The book covered many topics including: the mathematics of election, the mathematics of sharing, the mathematics of getting around, the mathematics of touring, the mathematics of networks, the mathematics of scheduling, financial mathematics, the mathematics of symmetry, and fractal geometry. The textbook emphasized the real life application of the different concepts.

B. *Treatment*

A thirty minute presentation about concept maps and how to construct them took place in the eleventh week, in which several examples of concept maps were presented. An additional 20 minutes were given to students to come up with a number of terms for the Euler Circuit. The students came up with the following terms: vertex, edge, adjacent edges, adjacent vertices, degree of vertices, connected graph, components, path, Circuit, Euler path, Euler Circuit and Famous examples such as Königsberg bridges puzzle.

These terms were written on the board and then students were asked to draw a concept map for the Euler Circuit using the terms that they came up with as well as any others they could think of. The final map that was drawn on the board was the result of the 20-minute discussion between the researcher and students. The discussion emphasized the connections between all the terms that the students mentioned and used in constructing their concept maps.

C. *Procedure*

During the first week of the semester, the researchers explained to the instructors the purpose of the study. Data were collected from a pre-treatment test that was conducted during the eleventh week of the semester, student at this point already covered the Euler Circuit concept, as well as from a post-treatment test during the fourteenth week of the semester. The pre-treatment test consisted of one question that aimed at collecting as much information as possible about the students' understanding of the Euler Circuit. The students were asked to answer the following question:

- What does the Euler Circuit mean to you? Please mention as many points as you can.

The same question was used for the post-treatment test which was given during the fourteenth week of the semester. Both tests were administered by one of the researchers and taken by students without any interference from the researcher.

VI. ANALYSIS & RESULTS

The collected data were organized in four categories based on students' responses in both the pre-treatment and post-treatment tests. The four different categories are: formal definition (D); the degree of vertices theorem (T); graphical meaning (G) and unclear statement (U).

The unclear statement category included the incorrect expressions that students used to describe the Euler Circuits and the expressions that did not fall under any of the other categories. Following are examples of direct quotes from students' answers in that category:

- "The graph is path starting and ending at the same point"
- "It means: the graph is connected"
- "The graph has vertices and edges which are adjacent"
- "You can go through each vertex one time only"

The objective of the data analysis was to gain as much information as possible about students' understanding of the Euler Circuit and to compare between each student's understanding before and after the treatment. Tables I and II give a summary of students' analysed responses.

TABLE I

THE NUMBER OF STUDENTS IN EACH CATEGORY BEFORE AND AFTER TREATMENT.

| Category | Before | After |
|---|---|---|
| D | 16 | 60 |
| T | 25 | 55 |
| G | 21 | 44 |
| U | 38 | 6 |
| Total | 100 | 165 |

As evident in Table I, more students emphasized the following concept categories in the post-treatment test compared to the pre-treatment test according to the ratios indicated between parentheses: D (60 to 16), T (55 to 25), G (44 to 21). The table also shows that before treatment the number of students who used unclear statements was 38 compared to 6 after treatment. This is an indication that after treatment, more students emphasized more categories of the Euler Circuit.

Table II represents the number of categories for students' responses before and after treatment. Before the treatment, only 15 students mentioned more than one representation for the Euler Circuit, the majority of them (69 out of 83) mentioned only one representation or none for the Euler Circuit. After the treatment, using concept maps, the majority of students (55 out of 83) mentioned two or more representations for the Euler Circuit.

TABLE II

NUMBER OF STUDENTS ACCORDING TO THE TOTAL NUMBER OF CATEGORIS

| Total number of Categories used | Before | After |
|---|---|---|
| 0 | 39 | 6 |
| 1 | 30 | 21 |
| 2 | 14 | 27 |
| 3 | 1 | 28 |
| Total | 83 | 83 |

VII. DISCUSSION

This study aimed at investigating the effect of using concept maps on advancing students' understanding of the Euler Circuit concept. The 83 participants took a pre-treatment test during the eleventh week of the semester and completed a test after the treatment during the fourteenth week. The treatment consisted of a training on the use of concept maps and coming up with concept maps to the Euler Circuit.

The images that students held about the Euler Circuit varied from one student to another, before and after treatment. As shown in Table I, before treatment, the majority of students (38 out of 83) had unclear images of the Euler Circuit. In addition to that, 21 out of 83 students had graphical images for the Euler Circuit. This indicates that students' knowledge was mostly limited to procedural knowledge. The formal definition came fourth in students' images (16 out of 83). In addition to that, 25 students mentioned the theorem that can be used to find the Euler Circuit.

After treatment, the most dominant image of the Euler Circuit was the formal definition (60 out of 83), which is a significant change compared to only 16 before treatment. Geometrical meaning of Euler Circuit came third in students' images (44 out of 83). Fifty five mentioned the theorem that can be used to find the Euler Circuit, and 6 students used unclear statements.

The results, as indicated in Table II, show that the majority of students after treatment (55 students) used multiple representations of the Euler Circuit compared to only 15 students before treatment (see Table II). Based on the data in Table II, the mean of the number of categories used by students before the treatment is 0.73 compared to the mean 1.92 after the treatment. This means that the number of categories used by the students after the treatment increased. This would indicate that there is a noticeable improvement in students' conceptual understanding after treatment in comparison with their conceptual understanding prior to the treatment. Hence the students' images of the Euler Circuit were richer after the use of concept maps. The findings of this study would suggest that the use of concept maps helped students develop and construct a more appropriate understanding of the Euler Circuit after using concept maps.